\documentclass[12pt]{article}
\usepackage{amsfonts}
\usepackage{mathrsfs}
\usepackage{amsmath,amsthm,amscd,amsfonts,mathrsfs,amssymb,graphicx,enumerate}
\usepackage{xypic}
\usepackage{makeidx}

\newtheorem{thm}{Theorem}[section]
\newtheorem{cor}[thm]{Corollary}
\newtheorem{lem}[thm]{Lemma}
\newtheorem{prop}[thm]{Proposition}
\theoremstyle{definition}
\newtheorem{defn}[thm]{Definition}

 \DeclareMathOperator{\Spec}{Spec}

\newcommand{\Z}{\ensuremath\mathbb{Z}}

\newcommand{\PP}{\ensuremath\mathbb{P}}
\newcommand{\calO}{\ensuremath\mathcal{O}}

\newcommand{\calF}{\ensuremath\mathscr{F}}

\newcommand{\Set}[2]{\left\{#1:#2\right\}}

\begin{document}
\title{Foliations and Rational Connectedness in Positive Characteristic}
\author{Mingmin Shen}
\date{}
\maketitle

\begin{abstract}
In this paper, the technique of foliations in characteristic $p$ is
used to investigate the difference between rational connectedness
and separable rational connectedness in positive characteristic. The
notion of being freely rationally connected is defined; a variety is
freely rationally connected if a general pair of points can be
connected by a free rational curve. It is proved that a freely
rationally connected variety admits a finite purely inseparable
morphism to a separably rationally connected variety. As an
application, a generalized Graber-Harris-Starr type theorem in
positive characteristic is proved; namely, if a family of varieties
over a smooth curve has the property that its geometric generic
fiber is normal and freely rationally connected, then it has a
rational section after some Frobenius twisting. We also show that a
freely rationally connected variety is simply connected.
\end{abstract}
\textbf{Key words:} Foliations, rationally connected, formal neighborhoods.\\
\textbf{Subject Class:} \textbf{14B20, 14G15}

\section{Introduction}
Since the work of J. Koll\'ar, Y. Miyaoka and S. Mori on rationally
connected varieties, it has been widely accepted that varieties with
lots of rational curves on them provide a good generalization of
projective spaces. We fix $k$ to be an algebraically closed field.
In this paper, a variety over $k$ is an integral, separated, finite
type $k$-scheme. We recall the following definitions, c.f.\cite{js},
\cite{Ko}.
\begin{defn}
Let $X$ be a quasi-projective variety over $k$. A \textit{rational curve}
on $X$ is a nonconstant
morphism $\phi:\PP^1\rightarrow X$. Let
$T_X=\mathscr{H}om_{\mathcal{O}_X}(\Omega_{X/k},\mathcal{O}_X)$ be
the tangent sheaf of $X$. If the image of $\phi$ is contained in the smooth locus
$X^{sm}$ of $X$, then we can split $\phi^*T_X$ into direct sum of line bundles
$$\phi^*T_X\cong
\sum_{i=1}^{n}\calO_{\PP^1}(a_i),
$$
where $n=\dim X$. We say that $\phi$ is \textit{free} if
$\text{Im}(\phi)\subset X^{sm}$ and $a_i\geq0$, for all $1\leq i\leq
n$. We say that $\phi$ is \textit{very free} if
$\text{Im}(\phi)\subset X^{sm}$ and $a_i\geq1$, for all $1\leq i\leq
n$. We say that $X$ is \textit{separably uniruled} if it has a free
rational curve and \textit{separably rationally connected} (SRC) if
it has a very free rational curve. We say that $X$ is
\textit{rationally connected} (RC) if there is a variety $Y$ and a
morphism $u:\PP^1\times Y\rightarrow X$ such that
$u^{(2)}:\PP^1\times\PP^1 \times Y\rightarrow X\times X$ is
dominant. In this case we also say that a general pair of points on
$X$ are connected by a rational curve.
\end{defn}

In addition to these, we make the following

\begin{defn}
Let $X$ be as above. We say that $X$ is \textit{freely rationally
connected} (FRC), if there is a variety $Y$ and a morphism
$u:\PP^1\times Y\rightarrow X$ such that $u^{(2)}:\PP^1\times\PP^1
\times Y\rightarrow X\times X$ is dominant and in addition, each
rational curve parametrized by $Y$ is a free rational curve on $X$.
\end{defn}

In general, SRC implies FRC and FRC implies RC. It turns out that
over an algebraically closed field of characteristic 0, a smooth
projective variety is RC if and only if it is SRC. However, these
notions are not equivalent in positive characteristic. Namely, in
\cite{Ko} V.5.19, there is an example of a smooth variety in
characteristic $p$, which is FRC but not SRC. In \cite{ghs}, T.
Graber, J. Harris and J. Starr prove that a rationally connected
fibration over a curve admits a section in characteristic 0. In
\cite{js}, A.J. de Jong and J. Starr generalize this result to
characteristic $p$ and the price they pay is that they have to
assume that a general fiber is SRC.

In this paper, we use the technique of foliations in positive characteristic
to investigate the gap between the above
notions. Our main theorem is the following.

\begin{thm}
Let $k$ be an algebraically closed field of characteristic $p$. Let $X/k$ be a
quasi-projective algebraic variety. Assume that $X$ is freely
rationally connected. Then there exists a separably rationally
connected variety $Y$ and a finite purely inseparable morphism
$f:X\rightarrow Y$. If $X$ is regular in codimension 1 (or normal), then so is
$Y$.
\end{thm}

We sketch the main idea of the proof. First we use the free rational
curves on $X$ to construct a canonical subsheaf $\mathscr{D}$ of
$T_X$. Then we prove that $\mathscr{D}$ is closed under the Lie
bracket and taking $p^{\text{th}}$ power if the characteristic is
$p$. This will be done in Section 2. This result holds in arbitrary
characteristic and only relies on $X$ being separably uniruled. In
Section 3, we construct the quotient of $X$ by $\mathscr{D}$ and
prove that the quotient is separably uniruled (resp. FRC) if $X$ is
so. As a result, if the quotient is not SRC, then we repeat the
construction above. In general, it is not guaranteed that this
procedure terminates with an SRC variety. In section 4, we
investigate the relation between the above procedure and the formal
neighborhood of a free rational curve. In particular, we show that
if the repeated construction does not terminate with an SRC variety,
then the global regular formal functions of the formal neighborhood
of a free rational curve form a power series ring. In section 5, we
prove the main theorem by showing the fact that the FRC condition
will force the quotient procedure to terminate with an SRC variety.
Then we give an application of the main theorem to prove the
following Graber-Harris-Starr type theorem.
\begin{thm}
Let $\pi: \mathscr{X}\to B$ be a proper flat family over a smooth
curve $B$, here everything is over an algebraically closed field $k$
of characteristic $p$. Assume that the geometric generic fiber of
$\mathscr{X}\to B$ is normal and freely rationally connected. Then
there is a morphism $s:B\to \mathscr{X}$ such that $\pi\circ
s=F_{\text{abs},B}^d$ for some $d\geq 0$, where
$F_{\text{abs},B}:B\to B$ is the absolute Frobenius morphism.
\end{thm}
An interesting consequence of the above theorem is
\begin{cor}
Let $X/k$ be a proper normal FRC variety over an algebraically
closed field $k$ of characteristic $p$. Then $X$ is simply
connected. Namely, the algebraic fundamental group $\pi_1(X)$ is
trivial.
\end{cor}

\textbf{Acknowledgement}: The author would like to thank his
advisor, Aise Johan de Jong, without whose careful and patient
instruction, the author would not have been able to carry out this
research. Many thanks to the referee who pointed out that FRC should
imply simply connectedness.

\section{The Foliation}

In this section we will fix the following notations. Let $X/k$ be a
quasi-projective algebraic variety of dimension $n$ over $k$ (of
arbitrary characteristic, but we will focus on positive
characteristic later). Assume that $X$ is separably uniruled. Then
there is a free rational curve $\phi:\PP^1\rightarrow X$ such that
$$
\phi^*T_X \cong\sum^r_{i=1}\calO_{\PP^1}(a_i)\bigoplus\calO_{\PP^1}^{n-r},
$$
where $a_i>0$ for all $1\leq i\leq r$.

\begin{defn}
The distinguished subsheaf $\sum^r_{i=1}\calO_{\PP^1}(a_i)$ of
$\phi^*T_X$ will be denoted by $\text{Pos}(\phi^*T_X)$. The
\textit{positive rank} of $X$ is defined to be the largest number
$r$ such that $r=\text{rank}(\text{Pos}(\phi^*T_X))$ for some free
rational curve $\phi$ on $X$. A free rational curve $\phi$ is said
to be \textit{maximally free} if Pos($\phi^*T_X$) has rank equal to
the positive rank $r$ of $X$. We use
$\text{Hom}^\text{m.free}(\PP^1,X)$ to denote the open subscheme of
$\text{Hom}^\text{free}(\PP^1,X)$ that parametrizes all maximally
free rational curves on $X$. Let $U_{m}$ be the open subvariety of
$X$ defined by the image of
$\PP^1\times\text{Hom}^\text{m.free}(\PP^1,X)\to X$.
\end{defn}

\begin{prop}
Let $X$ be as above, then for each closed point $x:
\Spec(k)\rightarrow U_m$, there is a well defined subspace
$\mathscr{D}(x)$ of $x^*T_X=T_{X,x}\otimes k(x)$ such that for every
maximally free rational curve $\phi:\PP^1\rightarrow X$ through
which $x$ factors as
$$
\xymatrix{\Spec(k) \ar[r]^i &\PP^1\ar[r]^\phi & X}
$$
we always have $i^*\text{Pos}(\phi^*T_X)=\mathscr{D}(x)\subset x^*T_X$.
\end{prop}

\textbf{Proof.} Let $\phi_1$ and $\phi_2$ be two maximally free
rational curves with $\phi_1(0)=\phi_2(0)=x$. So we have the
following diagram
$$
\xymatrix{
                      &\Spec(k)\ar[dl]_{i_1} \ar[d]^{x}\ar[dr]^{i_2}   &  \\
C_1\ar[r]_{\phi_1}   &X           &C_2 \ar[l]^{\phi_2} }
$$
where $C_1$ and $C_2$ are two $\PP^1$'s, with $i_1$ and $i_2$ being
the inclusions of the two origins. We need to show that
$i_1^*\text{Pos}(\phi_1^*T_X)=i_2^*\text{Pos}(\phi_2^*T_X)$ as
subspaces of $x^*T_X$. We glue $C_1$ and $C_2$ at the origins and
get a nodal curve $C$ and a morphism $\phi:C\rightarrow X$. We will
prove that the deformation of $C$ is unobstructed and a general
deformation of $C$ is a free rational curve $\varphi$. Assume that
$i_1^*\text{Pos}(\phi_1^*T_X)\neq i_2^*\text{Pos}(\phi_2^*T_X)$,
then we will show that Pos($\varphi^*T_X$) has rank $>r$ and this is
a contradiction.

Let $F:(\text{Sch}/k)^\circ\rightarrow \text{(Set)}$ be the
deformation functor defined in the following way: $F(S)$ consists of
isomorphism classes of diagrams
$$
\xymatrix{\mathscr{C}\ar[d]^\pi\ar[r]^f &X\\
             S &}
$$
where $\pi$ is a proper flat family of at worst nodal curves of
genus 0 over the scheme $S$. Hence we can view $\phi: C\rightarrow
X$ as an element in $F(\Spec(k))$. The cotangent complex of $\phi$,
$L_\phi^\ast$, is the complex $\{\xymatrix@C=0.5cm{
  0 \ar[r] & \phi^*\Omega_{X/k} \ar[r] & \Omega_{C/k} \ar[r] & 0
  }\}$ with the $\Omega_{C/k}$ term having degree 0. The first order
  deformation of $\phi$ is given by $\mathbb{E}\text{xt}^1_{\calO_C}(L_\phi^\ast,
  \calO_C)$ and the obstruction space lives in $\mathbb{E}\text{xt}^2_{\calO_C}(L_\phi^\ast,
  \calO_C)$, \cite{cox} \cite{lti}. The spectral sequence for hyper-extension
  groups gives the following long exact sequence

\begin{equation}
\xymatrix{0\ar[r] & \mathbb{E}\text{xt}^0_{\calO_C}(L_\phi^\ast,
  \calO_C) \ar[r] & \text{Ext}^0_{\calO_C}(\Omega_{C/k}, \calO_C)
  \ar[r] & \text{Ext}^0_{\calO_C}(\phi^*\Omega_{X/k}, \calO_C)\\
           \ar[r] & \mathbb{E}\text{xt}^1_{\calO_C}(L_\phi^\ast,
  \calO_C) \ar[r]^\rho & \text{Ext}^1_{\calO_C}(\Omega_{C/k}, \calO_C)
  \ar[r] & \text{Ext}^1_{\calO_C}(\phi^*\Omega_{X/k}, \calO_C)\\
           \ar[r] & \mathbb{E}\text{xt}^2_{\calO_C}(L_\phi^\ast,
  \calO_C) \ar[r] & \text{Ext}^2_{\calO_C}(\Omega_{C/k}, \calO_C)
  \ar[r] & \text{Ext}^2_{\calO_C}(\phi^*\Omega_{X/k},
  \calO_C)=0}\label{spectral}
\end{equation}

\begin{lem}(\cite{dm}) The following are true:\\
(i) $\text{Ext}^2_{\calO_C}(\Omega_{C/k}, \calO_C)=0$, \\
(ii) $\text{Ext}^1_{\calO_C}(\phi^*\Omega_{X/k}, \calO_C)=
\text{H}^1(C,\phi^*T_X)=0$
\end{lem}

From the above lemma and the hyper-extension spectral sequence
\eqref{spectral}, we get that
$\mathbb{E}\text{xt}^2_{\calO_C}(L_\phi^\ast,\calO_C)=0$ and hence
the deformation of $\phi$ is unobstructed. Let
$$
\xymatrix{
  C\ar[r]^i\ar[d] &\mathscr{C}\ar[d]^\pi\ar[r]^f &X\\
  b_0\ar[r] &B &
}
$$
be a deformation of $\phi$ over a smooth curve $B$, i.e. $f\circ
i=\phi$.

\begin{lem}
We can choose ($\mathscr{C}/B,f$) such that\\
(i) a general fiber of $\pi$ is a smooth rational curve.\\
(ii) $f|_{\mathscr{C}_b}$ is a free rational curve for general $b\in
B(k)$.
\end{lem}
To prove (i), we first consider the deformation of $C$ as a
$k$-scheme; the first order deformation is
$\text{Ext}^1_{\calO_C}(\Omega_{C/k},\calO_C)$ and the obstruction
lives in $\text{Ext}^2_{\calO_C}(\Omega_{C/k},\calO_C)=0$. The
local-global spectral sequence gives
$$
\xymatrix{
    0\ar[r]
    &\text{H}^1(C,\mathscr{H}om(\Omega_C,\calO_C))\ar[r]
    &\text{Ext}^1_{\calO_C}(\Omega_C,\calO_C)\\
    \ar[r]
    &\text{H}^0(C,\mathscr{E}xt^1_{\calO_C}(\Omega_C,\calO_C))
    \ar[r]
    &\text{H}^2(C,\mathscr{H}om_{\calO_C}(\Omega_C,\calO_C))=0
}
$$
Note that
$\text{H}^0(C,\mathscr{E}xt^1_{\calO_C}(\Omega_C,\calO_C))$
classifies first order local deformations and
$\text{Ext}^1_{\calO_C}(\Omega_C,\calO_C)$ classifies all first
order (global) deformations of $C/k$. Hence all local deformations
of $C/k$ come from global deformations. By embedding $C$ into
$\PP^2$ and using the exact sequence
\begin{equation}
\xymatrix{ 0\ar[r] &\mathscr{I}/\mathscr{I}^2 \ar[r]
&\Omega_{\PP^2}|_C \ar[r] &\Omega_C \ar[r] &0},\label{c}
\end{equation}
to do local computation, we find that
$\mathscr{E}xt^1_{\calO_C}(\Omega_C,\calO_C)=k(P)$ where $P$ is the
node of $C$ and it follows that
$\text{H}^0(C,\mathscr{E}xt^1_{\calO_C}(\Omega_C,\calO_C))= k$.
Actually, one can compute all the formal local deformations
explicitly. Let $R=k[[X,Y]]/(XY)$, all deformations of $R/k$ to
$k[\epsilon]$ are given by
$$\Set{k[\epsilon][[X,Y]]/(XY-\lambda\epsilon)}{\lambda\in k}.$$
This means that there exist local deformations that smooth out the
node and hence there exist deformations of $C/k$ that smooth out the
node. In the sequence \eqref{spectral}, the surjectivity of $\rho$
says that all deformations of $C/k$ comes from deformations of
$\phi:C/k\rightarrow X$ by forgetting the morphism to $X$. Thus
there is a deformation, ($\mathscr{C}/B,f$), of $\phi:C/k\rightarrow
X$ smoothing out the node of $C$. This proves (i) of the Lemma.

To prove (ii), we consider the locally free sheaf $f^*T_X$ on
$\mathscr{C}$. First, by shrinking $B$ and replacing $B$ by a finite
\'etale cover, we may assume that $\mathscr{C}/B$ has a section
$\sigma$ passing through a smooth point of $C=\mathscr{C}_{b_0}$. So
$\sigma$ defines a horizontal divisor $D$ on $\mathscr{C}$. Let
$\mathscr{E}=f^*T_X(-D)$ and we have the following short exact
sequence
$$
\xymatrix{ 0\ar[r] &\mathscr{E}|_C\ar[r]
&j_{1*}(\mathscr{E}|_{C_1})\oplus j_{2*}(\mathscr{E}|_{C_2})\ar[r]
&\mathscr{E}(P)\ar[r] &0},
$$
where $P$ is the nodal point of $C$ and $j_1:C_1\rightarrow C$ and
$j_2:C_2\rightarrow C$ are the inclusions. The associated long exact
sequence is
$$
\xymatrix{
  \text{H}^0(C_1,\mathscr{E}|_{C_1})\oplus\text{H}^0(C_2,\mathscr{E}|_{C_2})\ar[r]^\alpha
    &\mathscr{E}(P)=\mathscr{E}_P\otimes k(P)\ar[r]
    &\\
    \text{H}^1(C,\mathscr{E}|_{C})
  \ar[r]
    &\text{H}^1(C_1,\mathscr{E}|_{C_1})\oplus\text{H}^1(C_2,\mathscr{E}|_{C_2})=0
    &}
$$
where the last term is 0 since $C_1$ and $C_2$ are free and $\alpha$
is surjective since one of $\mathscr{E}|_{C_1}$ and
$\mathscr{E}|_{C_2}$ is globally generated. It follows that
$\text{H}^1(C,\mathscr{E}|_{C})=0$. By semicontinuity theorem
\cite{har} III.12, we get
$\text{H}^1(\mathscr{C}_b,\mathscr{E}|_{\mathscr{C}_b})=0$ for
general $b\in B(k)$. Since
$\mathscr{E}|_{\mathscr{C}_b}=f_b^*T_X(-1)$, we know that $f_b$ is
free and this proves (ii).

Now we are ready to prove that for a general deformation
$f_b:\mathscr{C}_b\rightarrow X$, the positive part of $f_b^*T_X$
has rank $>r$ under the assumption that
$i_1^*\text{Pos}(\phi_1^*T_X)\neq i_2^*\text{Pos}(\phi_2^*T_X)$. We
may assume that $\mathscr{C}/B$ has two horizontal divisors $D_1$ and $D_2$
corresponding to sections $\sigma_i$ passing through a point of $C_i$ not
equal to $P$. Let $\mathscr{E}=f^*T_X(-D_1-D_2)$,
$\calF=\mathscr{E}|_{\mathscr{C}_{b_0}}$, $\calF_1=\calF|_{C_1}$ and
$\calF_2=\calF|_{C_2}$, then we have the following short exact
sequence,
$$
\xymatrix{0\ar[r] &\calF\ar[r] &i_{1*}\calF_1\oplus
i_{2*}\calF_2\ar[r] &\calF(P)=\calF_P\otimes k(P)\ar[r] &0}.
$$
The induced long exact sequence is
$$
\xymatrix{0\ar[r]
  &\text{H}^0(C,\calF)\ar[r]
  &\text{H}^0(C_1,\calF_1)\oplus\text{H}^0(C_2,\calF_2)
  &\\
\ar[r]^\beta
  &T_X(P)\ar[r]
  &\text{H}^1(C,\calF)\ar[r]
  &0 }.
$$
Note that $\text{Im}(\beta)$ is the subspace spanned by
$\text{Pos}(\phi_1^*T_X)(P)$ and $\text{Pos}(\phi_2^*T_X)(P)$, which
has dimension $>r$ by assumption. It follows that $\dim
\text{H}^1(C,\calF)=n-\dim(\text{Im}(\beta))<n-r$. By
semicontinuity again, we have
$\text{H}^1(\mathscr{C}_{b},\mathscr{E}|_{\mathscr{C}_b}) <n-r$.
Since $\mathscr{E}|_{\mathscr{C}_b}=f_b^*T_X\otimes\calO(-2)$ by
construction, the positive part of $f_b^*T_X$ has rank greater
than $r$. $\square$\\

This proposition implies that we can define $\mathscr{D}(x)$, for
each point $x\in U_m(k)$, to be the subspace
$i^*\text{Pos}(\phi^*T_X)$ of $x^*T_X$, where $\phi:\PP^1\rightarrow
X$ is a maximally free rational curve with $\phi(0)=x$, and
$i:\Spec(k)\rightarrow\PP^1$ is the inclusion of the origin.
$$
\xymatrix{ &\PP^1\ar[dr]^\phi &\\
       \Spec(k)\ar[ur]^i\ar[rr]^x&&X }
$$
Next, we will use faithfully flat descent technique to glue all the
$\mathscr{D}(x)$ together to get a sub-bundle of $T_X$ on some open
set.

\begin{prop}
There is a nonempty open subset $\tilde{U}\subset U_m$ such that
$\{\mathscr{D}(x):x\in \tilde{U}(k)\}$ glue together to give a
subbundle $\mathscr{D}$ of $T_X|_{\tilde{U}}$, i.e.
$x^*\mathscr{D}=\mathscr{D}(x)$, for all $x\in \tilde{U}(k)$.
\end{prop}

\textbf{Proof.} It is known that the universal morphism
$F:\PP^1\times\text{Hom}^\text{m.free}(\PP^1,X)\to X$ is smooth,
\cite{Ko} II.3. Take a connected open subset
$W\subset\text{Hom}^\text{m.free}(\PP^1,X)$. Consider the smooth
morphism $\pi:\PP^1\times W\rightarrow X$, in particular, $\pi$ is
flat. By shrinking $W$, we may assume that $\pi^*T_X$ splits
uniformly as $\pi^*T_X=\sum_{i=0}^r\calO(a_i)\oplus\calO^{n-r}$
where $\calO=\calO_{\PP^1_W}$ and $a_i>0$. Denote the distinguished
sub-bundle $\sum_{i=1}^r\calO(a_i)$ of $\pi^*T_X$ by $\mathscr{V}$.
Let $\tilde{U}\subset X$ be the image of $\pi$. By construction, we
have $\tilde{U}\subset U_m$. Then $\pi:\PP^1\times W\rightarrow
\tilde{U}$ is faithfully flat. Next, we want to construct descent
data on $\mathscr{V}$ from that of $\pi^*T_X$. To do this, we denote
$\PP^1\times W$ by $Z$; $Z\times_\pi Z$ by $Z^{(2)}$ with
projections $p_1$ and $p_2$; $Z\times_\pi Z\times_\pi Z$ by
$Z^{(3)}$ with projections $p_{ij}$, $1\leq i< j\leq3$. Let
$\pi_1=\pi:Z\to \tilde{U}$, $\pi_2:Z^{(2)}\to \tilde{U}$ and
$\pi_3:Z^{(3)}\to \tilde{U}$ be the obvious morphisms. First, we
want to construct an isomorphism $\varphi:p_1^*\mathscr{V}\to
p_2^*\mathscr{V}$ in the following way. Note that both
$p_1^*\mathscr{V}$ and $p_2^*\mathscr{V}$ can be viewed as
sub-bundles of $\pi_2^*T_X$ that are local direct summands and by
construction of $\mathscr{V}$ and $\{\mathscr{D}(x)\}$, we have
\begin{equation}
p_1^*\mathscr{V}\otimes k(z) =p_2^*\mathscr{V}\otimes k(z),
\quad\text{for all closed points } z\in \tilde{Z}^{(2)}\label{red}
\end{equation}
as subspaces of $\pi_2^*T_X\otimes k(z)$. Since $Z^{(2)}$ is smooth
and hence reduced, the condition \eqref{red} implies that
$p_1^*\mathscr{V}=p_2^*\mathscr{V}$ as sub-bundles of $\pi_2^*T_X$
by Hilbert Nullstellensatz(\cite{har}, I.1). We define $\varphi$ to
be the identification. Our next step is to show that this
isomorphism $\varphi$ satisfies the cocycle condition, c.f.
\cite{blr} Chapter 6,
\begin{equation}
p_{13}^*\varphi =p_{23}^*\varphi\circ p_{12}^*\varphi
\end{equation}
To be more precise, we need to show that the following diagram is commutative
$$
\xymatrix{
    p_{12}^*p_1^*\mathscr{V}\ar[r]^{p_{12}^*\varphi}\ar[d]^{id}
       &p_{12}^*p_2^*\mathscr{V}\ar[r]^{id}
       &p_{23}^*p_1^*\mathscr{V}\ar[r]^{p_{23}^*\varphi}
       &p_{23}^*p_2^*\mathscr{V}\ar[d]^{id}\\
    p_{13}^*p_1^*\mathscr{V}\ar[rrr]^{p_{13}^*\varphi}
       & & &p_{13}^*p_2^*\mathscr{V}
}
$$
All the above terms, as sub-bundles of $\pi_3^*T_X$, are the same.
Note that $p_{13}^*\varphi-p_{23}^*\varphi\circ p_{12}^*\varphi$, as
a morphism from $p_{13}^*p_1^*\mathscr{V}$ to
$p_{13}^*p_2^*\mathscr{V}$, reduces to 0 on fibers and hence is 0 by
Hilbert Nullstellensatz. This shows the cocycle condition. By
construction, the inclusion $\mathscr{V}\hookrightarrow \pi^*T_X$ is
compatible with the descent data constructed on $\mathscr{V}$ and
the canonical one on $\pi^*T_X$. Hence by faithfully flat descent we
conclude that there is a sub-bundle $\mathscr{D}$ of
$T_X|_{\tilde{U}}$ with $\mathscr{V}=\pi^*\mathscr{D}\hookrightarrow
\pi^*T_X$. $\square$

\textbf{Remark:} In general, let $X$ be a noetherian scheme and
$\mathscr{E}$ be a coherent $\calO_X$-module. Let $\mathscr{F}$ be a
submodule of $\mathscr{E}|_U$ where $U\subset X$ is a dense open
subscheme. Then it is a standard fact that there is a maximal
submodule $\tilde{\mathscr{F}}$ of $\mathscr{E}$ which extends
$\mathscr{F}$. Actually, we can define $\tilde{\mathscr{F}}$ by $
\Gamma(V,\tilde{\mathscr{F}})=\Set{s\in \Gamma(V,\mathscr{E})}
{s|_{V\cap U}\in \Gamma(V\cap U,\mathscr{F})} $. In addition, if $X$
is integral and $\mathscr{F}\subset\mathscr{E}|_U$ is saturated,
then $\tilde{\mathscr{F}}$ is a saturated submodule of $\mathscr{E}$
that extends $\mathscr{F}$. Hence, in our situation, there is a
canonically defined saturated subsheaf, still denoted by
$\mathscr{D}$, of $T_X$ whose restriction to $\tilde{U}$ is the
subbundle constructed in the above proposition. From now on, we will
always use $\mathscr{D}$ or
$\mathscr{D}_X$ to denote this canonical subsheaf of $T_X$.\\

\begin{prop}
The subsheaf $\mathscr{D}$ satisfies the following two properties.\\
(i) $[\mathscr{D},\mathscr{D}] \subseteq \mathscr{D}$, i.e.
$\mathscr{D}$ is closed under the Lie Bracket.\\
(ii) $\mathscr{D}^p\subseteq\mathscr{D}$, i.e. $\mathscr{D}$ is
closed under taking $p^{\text{th}}$ power if $k$ is of
characteristic $p$.
\end{prop}
\textbf{Proof.} We use the same notions as in the previous proof.
Consider the following diagram
$$
\xymatrix{
\mathscr{D}\otimes_k\mathscr{D}\ar[r]^{[-,-]}\ar[dr]
    &T_X\ar[r]
    &T_X/\mathscr{D}\\
  &\mathscr{D}\otimes_{\calO_X}\mathscr{D}\ar[ur]_\rho &
},
$$
where $\rho$ is $\calO_X$-linear. Pullback $\rho$ by the faithfully
flat map $\pi$ we get
$$
\pi^*\rho:\pi^*(\mathscr{D}\otimes_{\calO_X}\mathscr{D})=\bigoplus_{i,j}\calO(a_i+a_j)\rightarrow
\pi^*(T_X/\mathscr{D})=\calO^{n-r}
$$
Since $a_i>0$, we know that $\pi^*\rho=0$ and hence
$\rho|_{\tilde{U}}=0$. Since $T_X/\mathscr{D}$ is torsion free, we
get $\rho=0$. This proves (i). Now assume that $k$ is of
characteristic $p$ and consider the following composition $\varrho$
$$
\varrho:\xymatrix{\mathscr{D}\ar[r]^{(-)^p} &T_X\ar[r] &T_X/\mathscr{D}}.
$$
It is known that $\varrho$ is $p$-linear, i.e.
$\varrho(fD)=f^p\varrho(D)$ where $f\in\calO_X$, if $\mathscr{D}$ is
closed under Lie bracket, \cite{kat}. This defines an
$\calO_X$-linear morphism $\varrho'=1\otimes\varrho:
F_{\text{abs},X}^*\mathscr{D} \to T_X/\mathscr{D}$, where
$F_{\text{abs},X}$ is the absolute Frobenius morphism of $X$. If we
pull the diagram back via $\pi$ and note that
$\pi^*F_{\text{abs},X}^*\mathscr{D}=F_{\text{abs},Z}^*\pi^*\mathscr{D}$
and $F_{\text{abs},Z}^*\calO(a) =\calO(p a)$, we get
$$
\pi^*\varrho':\pi^*(F_{\text{abs},X}^*\mathscr{D})=\bigoplus_{i=1}^r
\calO(p a_i) \to \pi^*(T_X/\mathscr{D})=\calO^{n-r}.
$$
Then $\pi^*\varrho'$ has to be 0 and hence $\varrho'=0$. This proves
(ii). $\square$

\section{The Quotient}
In this section, we assume that $k$ is of characteristic $p$. Let
$X/k$ be a separably uniruled variety. Let $r$ be the positive rank
of $X$. We apply the results from the previous section and define
the foliation $\mathscr{D}$ on $X$. We define a sheaf of
$k$-algebras $\mathscr{A}\hookrightarrow \calO_X$ by
$$
\Gamma(V,\mathscr{A})=\{f\in\Gamma(V,\calO_X)|Df=0, \forall
D\in\Gamma(V_0,\mathscr{D}),\forall\text{ open } V_0\subset V  \}
$$
Note that we have the following inclusions
$$
\xymatrix{
  \calO_X\ar[r]^{(-)^p} &\mathscr{A}\ar[r] &\calO_X
}.
$$
In this way, $\mathscr{A}$ can be viewed as a sheaf of
$\calO_X$-algebras. We define
$Y:=\textbf{Spec}_{\calO_X}(\mathscr{A})$. If we define the relative
Frobenius morphism $F_{X/k}:X\to X^{(1)}$ in the following way,
$$
\xymatrix{
X\ar[drr]^{F_{\text{abs},X}}\ar[ddr]\ar[dr]
   &
   &\\
   &X^{(1)}\ar[r]^\sigma\ar[d]
   &X\ar[d]\\
   &\Spec(k)\ar[r]^{F_{\text{abs},k}}
   &\Spec(k)
}
$$
then we have the natural morphisms
\begin{equation}
\xymatrix{X\ar[r]^f &Y\ar[r]^g &X^{(1)}\ar[r]^\sigma &X}
\end{equation}
where $f$ and $g$ are $k$-linear and $g\circ f=F_{X/k}$. We need the
following

\begin{prop}[\cite{eke}, \cite{miy}]
Let $X$ and $Y$ be as above, then\\
(i) On any open part $U$ where $\mathscr{D}$ is a subbundle, $Y$ is
smooth; for each closed point $x\in U$, there are formal coordinates
$t_1,t_2,\ldots, t_n$ such that formally locally $f$ is given by
$(t_1,\ldots,t_r,t_{r+1},\ldots,t_n)\mapsto
(t_1^p,\ldots,t_r^p,t_{r+1},\ldots,t_n)$ and $\hat{\mathscr{D}}_x$
is freely generated by $\{\frac{\partial}{\partial t_1},\ldots,
\frac{\partial}{\partial t_r}\}$. In particular, $f$
is faithfully flat on $U$.\\
(ii) There is a canonical exact sequence
\begin{equation}
\xymatrix{
  0\ar[r] &\mathscr{D}\ar[r] &T_X\ar[r]^{df} &f^*T_Y\ar[r]^{f^*\alpha} &F_{\text{abs},X}^*\mathscr{D}
\ar[r] &0
}\label{fuda}
\end{equation}
where $\alpha$ is induced by $dg$. \\
(iii) Let $X$ be a normal variety. Then there is a one-to-one
correspondence between the foliations $\mathscr{D}$ on $X$ and the
normal varieties $Y$ between $X$ and $X^{(1)}$. \label{ekedahl}
\end{prop}

\begin{cor}\label{coreke}
If $X$ is regular in codimension 1, then so is $Y$. If $X$ is
normal, then so is $Y$.
\end{cor}

Now let $\phi:\PP^1\to X$ be a maximally free rational curve on $X$.
We have $d\phi:T_{\PP^1}\rightarrow \phi^*T_X$ is nonzero and
$T_{\PP^1}\cong \calO(2)$, hence
$d\phi(T_{\PP^1})\subset\text{Pos}(\phi^*T_X)$. Thus
$d(f\circ\phi)=0$, and hence $f\circ\phi$ factors through the
relative Frobenius of $\PP^1$, i.e.
$$
\xymatrix{
  X\ar[r]^f
     &Y\ar[r]^g
     &X^{(1)}\ar[r]^\sigma
     &X \\
  \PP^1\ar[u]^\phi\ar[r]_{F_{\PP^1/k}}
     &\PP^1\ar[u]_{\tilde{\phi}}\ar[r]^\sigma
     &\PP^1\ar[ur]_\phi
     &
}
$$
If we pull back the exact sequence \eqref{fuda} to the rational
curves and note that $\phi^*f^*=F_{\PP^1/k}^*\tilde{\phi}^*$, we get
\begin{equation}
\xymatrix{ 0\ar[r] & \phi^*\mathscr{D}\ar[r] & \phi^*T_X\ar[r]
&F_{\PP^1/k}^*(\tilde{\phi}^*T_Y)\ar[r]
 & F_{\PP^1/k}^*(\tilde{\phi}^*g^*\sigma^*\mathscr{D})\ar[r] &0
}\label{plbk}
\end{equation}
Let $\mathscr{Q}=\ker(\alpha)$, i.e.
\begin{equation}
\xymatrix{0\ar[r] &\mathscr{Q}\ar[r] &T_Y\ar[r]^\alpha
&g^*\sigma^*\mathscr{D}\ar[r] &0}\label{splt}
\end{equation}
It is easy to see that $\phi^*\mathscr{D}=\bigoplus_{i=1}^r
\calO(a_i)$ with $a_i>0$, and $\phi^*f^*\mathscr{Q}=\calO^{n-r}$
(which implies that $\tilde{\phi}^*\mathscr{Q}=\calO^{n-r}$) and
that $\tilde{\phi}^*(g^*\sigma^*\mathscr{D})
=\sigma^*\phi^*\mathscr{D}=\bigoplus_{i=1}^r\calO(a_i)$. Then the
$\tilde{\phi}$ pullback of sequence \eqref{splt} becomes
\begin{equation}
\xymatrix{ 0\ar[r] & \calO^{n-r}\ar[rr]
&&\tilde{\phi}^*T_Y\ar[rr]^{\tilde{\phi}^*\alpha}
&&\bigoplus_{i=1}^r\calO(a_i)\ar[r] &0 }\label{eqn2}
\end{equation}
From the above exact sequence, it is very easy to see that
$\tilde{\phi}$ is a free rational curve on $Y$ and hence $Y$ is
separably uniruled.

\begin{lem}
In the exact sequence \eqref{eqn2}, we have either the rank of
Pos($\tilde{\phi}^*T_Y$) is greater than $r$ or the sequence is
splitting.\label{rank_arg}
\end{lem}
\textbf{Proof.} Consider the following diagram
$$
\xymatrix{
    &0   &0   &   &\\
    &\calO^{n-r'}\ar[u]\ar[r]^{id} &\calO^{n-r'}\ar[u] &   &\\
0\ar[r] &\calO^{n-r}\ar[u]\ar[r] &\tilde{\phi}^*T_Y\ar[u]\ar[r]
&\sum_{i=0}^{r}\calO(a_i)\ar[r]  &0\\
0\ar[r] &\mathscr{E}\ar[u]\ar[r]
&\text{Pos}(\tilde{\phi}^*T_Y)\ar[r]\ar[u]
&\sum_{i=0}^{r}\calO(a_i)\ar[u]_{id}\ar[r] &0\\
  & 0\ar[u] &0\ar[u] & &
 }
$$
It follows easily that $\mathscr{E}\cong \calO^{r'-r}$ and hence the
lemma. $\square$\\

\begin{prop}
Let $X$ and $Y$ be as above. Then the following are true\\
(a) We have exactly one of the following cases\\
\indent (``Trivial case") The variety $X$ is SRC and $Y=X^{(1)}$\\
\indent
(``General case") The positive rank of $Y$ is strictly
greater than the positive rank of $X$.\\
\indent (``Splitting case") The variety $X$ is not SRC and the
positive rank of $Y$ is equal to the positive rank of $X$. In this
case, the exact sequence \eqref{splt}
splits canonically on some nonempty open set $\tilde{U}$ that contains a maximally free rational curve.\\
(b) If $X$ is FRC then so is $Y$. \label{cases}
\end{prop}
\textbf{Proof.} For (a), we only need to prove that the exact
sequence \eqref{splt} splits canonically in the ``Splitting case".
The above Lemma \eqref{rank_arg} shows that in the ``Splitting case"
$\tilde{\phi}^*T_Y$ has the same splitting type as $\phi^*T_X$, for
all maximally free rational curves $\phi$. Then we know that the
pull-back via $\tilde{\phi}$ of the exact sequence \eqref{splt} has
a unique splitting which identifies
$\tilde{\phi}^*(g^*\sigma^*\mathscr{D})$ as the sub-bundle
$\tilde{\phi}^*\mathscr{D}_Y$ of $\tilde{\phi}^*T_Y$, where
$\mathscr{D}_Y$ is the foliation on $Y$ as constructed in the
previous section. We can find an open subset $W$ of
$\text{Hom}^\text{free}(\PP^1,X)$, such that $(f\circ F)^*T_Y$
splits uniformly on $\PP^1\times W$ where $F:\PP^1\times W\to X$ is
the natural morphism. The composition $\mathscr{D}_Y \hookrightarrow
T_Y\to g^*\sigma^*\mathscr{D}$ is an isomorphism after pulling back
by $f\circ F$, hence itself is an isomorphism by faithfully flat
descent. This proves that the exact sequence \eqref{splt} splits
canonically on $\tilde{U}$. To prove (b), we only need to use the
fact that $f:X\to Y$ is dominant. $\square$

Proposition \eqref{cases} enables us to repeat the construction of
foliation on the quotients and we will get a sequence, which will be
called the quotient sequence of $X$.
$$
\xymatrix{
X\ar[r]^f &Y_1\ar[r]^{f_1} &Y_2\ar[r]^{f_2} &Y_3\ar[r]^{f_3} &\cdots
}
$$
\begin{cor}
If the ``splitting case" does not happen for infinitely many times
in the above procedure, then there is some $N\geq 1$ such that $Y_i$
is SRC, $Y_{i+1}=Y_i^{(1)}$ and $f_i$ is the relative Frobenius
morphism, for all $i\geq N$.\label{stops}
\end{cor}
\textbf{Proof.} This is because the positive rank of $Y_i$ is
bounded above by $n=\dim(X)$. $\square$

\section{The Formal Neighborhood}
In this section, we assume that the characteristic of $k$ is $p$.
let $X/k$ be a separably uniruled variety. We have the quotient
sequence
\begin{equation}
\xymatrix{ X\ar[r]^f &Y_1\ar[r]^{f_1} &Y_2\ar[r]^{f_2}
&Y_3\ar[r]^{f_3} &\cdots }\label{quot}
\end{equation}
We want to relate the ``splitting case" to the formal neighborhood
of a free rational curve on $X$. Before doing that, we state some
basic facts that are true in any characteristic.

\begin{lem}
Let $i:Z\to Y$ be a closed immersion of a smooth variety $Z/k$ into
another smooth variety $Y/k$. Let $\mathscr{I}$ be the
ideal sheaf defining $Z$ in $Y$. Then\\
 (i) The sequence $\xymatrix{0\ar[r]
&\mathscr{I}/\mathscr{I}^2\ar[r] &i^*\Omega_{Y/k}\ar[r]
&\Omega_{Z/k}\ar[r] &0}$
is exact and locally splitting.\\
 (ii) The sheaf $\mathscr{I}/\mathscr{I}^2$ is a locally free sheaf
on $Z$ whose rank is $\dim(Y)-\dim(Z)$.\\
 (iii) As a sheaf on $Z$, $\mathscr{I}$ locally generated by a
regular sequence of length $\dim(Y)-\dim(Z)$; in particular, $Z$ is
a local complete
intersection in $Y$.\\
 (iv) The sheaf $i^*\mathscr{I}^d=\mathscr{I}^d/\mathscr{I}^{d+1}$
is canonically isomorphic to $Sym^d(\mathscr{I}/\mathscr{I}^2)$, the
$d^\text{th}$ symmetric power of $\mathscr{I}/\mathscr{I}^2$, for
$d\geq 1$.
\end{lem}
\textbf{Proof.} These facts are standard.\\

Now let $f:X\to Y=Y_1$ be the quotient by $\mathscr{D}$ and $\phi$
be a maximally free rational curve on $X$ with the induced one on
$Y$ being $\tilde{\phi}$ as before. Then we have
\begin{lem}
If $\phi:\PP^1\to X$ is a closed immersion, then so is
$\tilde{\phi}:\PP^1\to Y$.
\end{lem}
\textbf{Proof.} Consider the following diagram
\begin{equation*}
\xymatrix
{
   X \ar[r]^f &Y \ar[r]^g &X^{(1)}\\
   \PP^1\ar[r]^F \ar[u]^\phi &\PP^1\ar[u]^{\tilde{\phi}} \ar[ur]_{\phi^{(1)}}  &
}
\end{equation*}
By assumption, both $\phi$ and $\phi^{(1)}$ are closed immersions.
It is easy to check that $\tilde{\phi}$ is also a closed
immersion. $\square$\\

\textbf{Notation and assumptions:} Now consider the following
situation
\begin{equation}
\xymatrix{
  X\ar[r]^f
     &Y\ar[r]^g
     &X^{(1)}\ar[r]^\sigma
     &X \\
  \PP^1\ar[u]^\phi\ar[r]_{F_{\PP^1/k}}
     &\PP^1\ar[u]_{\tilde{\phi}}\ar[r]^\sigma
     &\PP^1\ar[ur]_\phi
     &
}\label{stp1}
\end{equation}
where $X\to Y$ is the quotient by $\mathscr{D}$ and $\phi$ is a
maximally free rational curve on $X$. Assume that $\phi$ is a closed
immersion, then so is $\tilde{\phi}$ by the above lemma. Let
$\mathscr{I}$ be the ideal sheaf defining $\phi$ and
$\tilde{\mathscr{I}}$ be the ideal sheaf defining $\tilde{\phi}$.
\begin{lem}
In the ``splitting case" we have $\mathscr{I}/\mathscr{I}^2\cong
\calO^{n-r}\oplus \mathscr{N}$ and
$\tilde{\mathscr{I}}/\tilde{\mathscr{I}}^2\cong \calO^{n-r}\oplus
\tilde{\mathscr{N}}$ where both $\mathscr{N}$ and
$\tilde{\mathscr{N}}$ are direct sums of line bundles of negative
degree.\label{conormal}
\end{lem}
\textbf{Proof.} This follows from the short exact sequence
$$\xymatrix{0\ar[r] &\mathscr{I}/\mathscr{I}^2\ar[r] &\phi^*\Omega_X\ar[r]
&\Omega_{\PP^1}\ar[r] &0}$$
 and the corresponding one for $\tilde{\phi}$. $\square$\\

Note that $\phi$ and $\tilde{\phi}$ have the same underlying
topological space, which will be denoted by $\PP^1$ by abuse of
notation. Since in the diagram \eqref{stp1} the relative Frobenius
morphism $F_{\PP^1/k}$ induces an injective ring homomorphism,
$F_{\PP^1/k}^*:\calO_{\PP^1}\to (F_{\PP^1/k})_*\calO_{\PP^1}$, on
the structure sheaves and the functions on $Y$ are exactly the
functions on $X$ that are killed by $\mathscr{D}$, it follows that
the sheaf $\tilde{\mathscr{I}}$ consists of elements of the sheaf
$\mathscr{I}$ that are killed by $\mathscr{D}$. Hence we have a
natural morphism of abelian sheaves $\theta:
\tilde{\mathscr{I}}/\tilde{\mathscr{I}}^2\rightarrow
\mathscr{I}/\mathscr{I}^2$. Note that
$F_{\PP^1/k}^*\tilde{\phi}^*\Omega_Y=\calO_{\PP^1}\otimes_{\calO_{\PP^1}}\tilde{\phi}^*\Omega_Y$
and we have a natural sheaf homomorphism $1\otimes
id:\tilde{\phi}^*\Omega_Y\to F_{\PP^1/k}^*\tilde{\phi}^*\Omega_Y$.
Let $\vartheta$ be the composition:
$$
\xymatrix{ \tilde{\phi}^*\Omega_Y \ar[rr]^{1\otimes id}
&&F_{\PP^1/k}^*\tilde{\phi}^*\Omega_Y \ar[rr]^{\phi^*(df)^{\vee}}
&&\phi^*\Omega_X }
$$

\begin{lem}
Notation as above, we have:\\
(i) The following diagram, as abelian sheaves on the underlying
topological space of $\PP^1$, is commutative.
\begin{equation}
\xymatrix{
0\ar[r] &\mathscr{I}/\mathscr{I}^2\ar[r] &\phi^*\Omega_X\ar[r]
      &\Omega_{\PP^1}\ar[r] &0\\
0\ar[r]
&\tilde{\mathscr{I}}/\tilde{\mathscr{I}}^2\ar[r]\ar[u]^{\theta}
&\tilde{\phi}^*\Omega_Y\ar[r]\ar[u]^{\vartheta}
      &\Omega_{\PP^1}\ar[r]\ar[u]^{0} &0
}\label{commutative}
\end{equation}
(ii) The map
$\text{H}^1(\theta):\text{H}^1(\PP^1,\tilde{\mathscr{I}}
/\tilde{\mathscr{I}}^2)\to
\text{H}^1(\PP^1,\mathscr{I}/\mathscr{I}^2)$ is 0.\\
(iii) Let
$\theta_d:\tilde{\mathscr{I}}^d/\tilde{\mathscr{I}}^{d+1}\to
\mathscr{I}^d/\mathscr{I}^{d+1}$ be the natural morphism, then
$\text{H}^1(\theta_d):\text{H}^1(\PP^1,\tilde{\mathscr{I}}^d/
\tilde{\mathscr{I}}^{d+1})\to
\text{H}^1(\PP^1,\mathscr{I}^d/\mathscr{I}^{d+1})$ is
0.\label{vanishing}
\end{lem}
\textbf{Proof.} (i) is direct checking of the definitions of the
sheaf homomorphisms involved. To prove (ii), we note that the
diagram factors as the following diagram
$$
\xymatrix{
0\ar[r] &\mathscr{I}/\mathscr{I}^2\ar[r] &\phi^*\Omega_X\ar[r]
      &\Omega_{\PP^1}\ar[r] &0\\
0\ar[r]
&F^*(\tilde{\mathscr{I}}/\tilde{\mathscr{I}}^2)\ar[r]\ar[u]^{\tilde{\theta}}
      &F^*\tilde{\phi}^*\Omega_Y\ar[r]\ar[u]_{\phi^*(df)^{\vee}}
      &F^*\Omega_{\PP^1}\ar[r]\ar[u]_{dF=0} &0\\
0\ar[r] &\tilde{\mathscr{I}}/\tilde{\mathscr{I}}^2\ar[r]\ar[u]
       &\tilde{\phi}^*\Omega_Y\ar[r]\ar[u]
      &\Omega_{\PP^1}\ar[r]\ar[u] &0
}
$$
where $F=F_{\PP^1/k}$. Hence, we only need to show that
$\text{H}^1(\tilde{\theta})=0$. To do this, we consider the
following diagram
\begin{equation}
\xymatrix{
   &   &0   &   &\\
   &   &\phi^*\mathscr{D}^\vee\ar[u] &  &\\
0\ar[r] &\mathscr{I}/\mathscr{I}^2\ar[r]^{\lambda}
&\phi^*\Omega_X\ar[r]\ar[u]
      &\Omega_{\PP^1}\ar[r] &0\\
0\ar[r]
&F^*(\tilde{\mathscr{I}}/\tilde{\mathscr{I}}^2)\ar[r]^{\mu}\ar[u]^{\tilde{\theta}}
      &F^*\tilde{\phi}^*\Omega_Y\ar[r]\ar[u]_{\phi^*(df)^{\vee}}\ar@{.>}[ul]_\tau
      &F^*\Omega_{\PP^1}\ar[r]\ar[u]_{dF=0} &0\\
   &   & F_{\text{abs}}^*\phi^*\mathscr{D}^\vee\ar[u] &  &\\
   &   &0\ar[u]   &   &
}\label{bigdiagram}
\end{equation}
where the existence of $\tau$ is guaranteed by the fact that $dF=0$
and the resulting diagram is commutative. It is easy to see from the
above diagram that $\tau$ factors through
$$
\text{coker}(F_{\text{abs}}^*\phi^*\mathscr{D}^\vee\rightarrow
F^*\tilde{\phi}^*\Omega_Y)=\phi^*f^*\mathscr{Q}^{\vee}=\calO^{n-r}
$$
This means that $\tau$ factors as $F^*\tilde{\phi}^*\Omega_Y\to
\calO^{n-r}\to \mathscr{I}/\mathscr{I}^2$. In particular, this
implies that $\tilde{\theta}$ factors as
$F^*(\tilde{\mathscr{I}}/\tilde{\mathscr{I}}^2)\to\calO^{n-r}\to
\mathscr{I}/\mathscr{I}^2$ and hence $\text{H}^1(\tilde{\theta})=0$.
For (iii), first we factorize $\theta_d$ as
$$
\xymatrix{ \tilde{\mathscr{I}}^d/\tilde{\mathscr{I}}^{d+1}\ar[r]
&F^*(\tilde{\mathscr{I}}^d/\tilde{\mathscr{I}}^{d+1})
\ar[rr]^{\tilde{\theta}_d} &&\mathscr{I}^d/\mathscr{I}^{d+1} }
$$
and we only need to show that $\text{H}^1(\tilde{\theta}_d)=0$. To
do this, we consider the following commutative diagram
$$
\xymatrix{ F^*(\tilde{\mathscr{I}}^d/\tilde{\mathscr{I}}^{d+1})
\ar[rr]^{\tilde{\theta}_d}
&&\mathscr{I}^d/\mathscr{I}^{d+1}\\
   Sym^d(F^*(\tilde{\mathscr{I}}/\tilde{\mathscr{I}}^2))\ar[rr]^{Sym^d(\tilde{\theta})}\ar[u]\ar[dr]
&&Sym^d(\mathscr{I}/\mathscr{I}^{2})\ar[u]\\
   &Sym^d(\calO^{n-r})\ar[ur] &
}
$$
Since the vertical arrows are isomorphisms and
$Sym^d(\tilde{\theta})$ factors through $Sym^d(\calO^{n-r})$, we get
$\text{H}^1(\tilde{\theta}_d)=0$. $\square$

\begin{lem}
Notation and assumptions as above. Let
$\calO_Y/\tilde{\mathscr{I}}^d\to\calO_X/\mathscr{I}^d$ be the
natural morphism of abelian sheaves on the underlying topological
space of $\PP^1$ which induces the maps
$\alpha_d:\text{H}^0(\PP^1,\calO_Y/\tilde{\mathscr{I}}^d) \to
\text{H}^0(\PP^1,\calO_X/\mathscr{I}^d)$. Then $\alpha_d(a)$ lifts
to $\text{H}^0(\PP^1,\calO_X/\mathscr{I}^{d+1})$ for all $a\in
\text{H}^0(\PP^1,\calO_Y/\tilde{\mathscr{I}}^d)$.\label{lifting}
\end{lem}
\textbf{Proof.} This is a direct application of
Lemma\eqref{vanishing}. Consider the following short exact sequences
of abelian sheaves.
\begin{equation}
\xymatrix{
  0\ar[r] &\mathscr{I}^d/\mathscr{I}^{d+1}\ar[r] &\calO_X/\mathscr{I}^{d+1}\ar[r] &\calO_X/\mathscr{I}^d\ar[r] &0\\
  0\ar[r] &\tilde{\mathscr{I}}^d/\tilde{\mathscr{I}}^{d+1}\ar[u]_{\theta_d}\ar[r]
           &\calO_Y/\tilde{\mathscr{I}}^{d+1}\ar[r]\ar[u] &\calO_Y/\tilde{\mathscr{I}}^d\ar[r]\ar[u] &0
}\label{xysequence}
\end{equation}
The associated long exact sequences are
\begin{equation}
\xymatrix@C=0.5cm{
  0\ar[r]
    &\text{H}^0(\mathscr{I}^d/\mathscr{I}^{d+1})\ar[r]
    &\text{H}^0(\calO_X/\mathscr{I}^{d+1})\ar[r]
    &\text{H}^0(\calO_X/\mathscr{I}^d)\ar[r]^{\delta_d}
    &\text{H}^1(\mathscr{I}^d/\mathscr{I}^{d+1})\\
  0\ar[r]
    &\text{H}^0(\tilde{\mathscr{I}}^d/\tilde{\mathscr{I}}^{d+1})\ar[r]\ar[u]_{H^0(Sym^d(\theta))}
    &\text{H}^0(\calO_Y/\tilde{\mathscr{I}}^{d+1})\ar[r]\ar[u]_{\alpha_{d+1}}
    &\text{H}^0(\calO_Y/\tilde{\mathscr{I}}^d)\ar[r]^{\tilde{\delta}_d}\ar[u]_{\alpha_d}
    &\text{H}^1(\tilde{\mathscr{I}}^d/\tilde{\mathscr{I}}^{d+1})\ar[u]_{\text{H}^1(\theta_d)=0}
 }\label{xysequencel}
\end{equation}

The result is easy diagram chasing. $\square$

\begin{lem}
Notation as in the previous lemma and we further assume that $X\to
Y$ is of ``splitting case". If $\alpha_d$ is surjective (resp.
isomorphism) and $\tilde{\delta}_d=0$ in the diagram
\eqref{xysequencel} then $\alpha_{d+1}$ is also surjective (resp.
isomorphism) and $\delta_d=0$.\label{iso}
\end{lem}
\textbf{Proof.} First we claim that in the splitting case
$H^0(Sym^d(\theta))$ is an isomorphism for all $d\geq 1$. By Lemma
\eqref{conormal}, we have $\mathscr{I}/\mathscr{I}^2\cong
\calO^{n-r}\oplus \mathscr{N}$ and
$\tilde{\mathscr{I}}/\tilde{\mathscr{I}}^2\cong \calO^{n-r}\oplus
\tilde{\mathscr{N}}$ where $\mathscr{N}$ and $\tilde{\mathscr{N}}$
are direct sums of line bundles of negative degrees. From the
diagram \eqref{bigdiagram} we see that the sheaf homomorphism
$\theta$ maps $\calO^{n-r}\subset
\tilde{\mathscr{I}}/\tilde{\mathscr{I}}^2$ isomorphically to
$\calO^{n-r}\subset F_*(\calO^{n-r})\subset
F_*(\mathscr{I}/\mathscr{I}^2)$ and is zero on
$\tilde{\mathscr{N}}$. Then it is easy to see that
$\mathscr{I}^d/\mathscr{I}^{d+1}\cong
Sym^d(\mathscr{I}/\mathscr{I}^2)\cong Sym^d(\calO^{n-r})\oplus
\mathscr{N}_d$ and $\tilde{\mathscr{I}}^d/\tilde{\mathscr{I}}^{d+1}
\cong Sym^d(\tilde{\mathscr{I}}/\tilde{\mathscr{I}}^2)\cong
Sym^d(\calO^{n-r})\oplus \tilde{\mathscr{N}}_d$ where
$\mathscr{N}_d$ and $\tilde{\mathscr{N}}_d$ are direct sums of line
bundles of negative degrees. Hence the sheaf homomorphism
$Sym^d(\theta)$ factors as
$\tilde{\mathscr{I}}^d/\tilde{\mathscr{I}}^{d+1} \twoheadrightarrow
Sym^d(\calO^{n-r})\hookrightarrow F_*(Sym^d(\calO^{n-r}))\subset
F_*(\mathscr{I}^d/\mathscr{I}^{d+1})$. Note that here we identify a
sheaf $\calF$ on $\PP^1$ with $F_*(\calF)$. Then it is easy to see
that $H^0(Sym^d(\theta))$ is an isomorphism for all $d\geq 1$. Then
the Lemma is an easy diagram chasing in \eqref{xysequencel}.
$\square$\\

Now we are ready to prove the main theorem of this section.
\begin{thm}
Let $X/k$ be an $n$ dimensional uniruled algebraic variety over an
algebraically closed field $k$ of characteristic $p>0$. Assume that
the quotient sequence \eqref{quot} of $X$ has ``splitting case" in
each step. Let $\phi=\phi_0:\PP^1\to X$ be a maximally free rational
curve on $X$, which is a closed immersion, and $\phi_i$ be the
induced one on $Y_i$. Let $\mathfrak{X}=X_{/\PP^1}$ be the formal
neighborhood of $\phi$ in $X$ and $\mathfrak{Y}_i=Y_{i/\PP^1}$ be
the formal neighborhood
of $\phi_i$ in $Y_i$. Then\\
(a)$\Gamma(\mathfrak{Y}_{i+1},\calO_{\mathfrak{Y}_{i+1}})\to
\Gamma(\mathfrak{Y}_{i},\calO_{\mathfrak{Y}_{i}})$ is an isomorphism
for all $i\geq0$, where $\mathfrak{Y}_0=\mathfrak{X}$.\\
(b)$\Gamma(\mathfrak{X},\calO_{\mathfrak{X}})\cong
k[[t_1,\cdots,t_{n-r}]]$ is a formal power series ring of $n-r$
variables, where $r$ is the positive rank of $X$.\label{nbhd}
\end{thm}
\textbf{Proof.} Let $\mathscr{I}_i$ be the ideal sheaf on $Y_i$ that
defines $\phi_i$ as a closed subvariety of $Y_i$, for
$i=0,1,\cdots$, here $Y_0=X$ and $\phi_0=\phi$; we also write
$\mathscr{I}_0$ as $\mathscr{I}$. Our first claim is that
$\pi_{i,d}:\text{H}^0(\PP^1,\calO_{Y_i}/\mathscr{I}_i^{d+1})\rightarrow
\text{H}^0(\PP^1,\calO_{Y_{i-1}}/\mathscr{I}_{i-1}^{d+1})$ is
isomorphism for all $i\geq 1$, $d\geq0$. To prove this, we let
$\delta_{i,d}:\text{H}^0(\PP^1,\calO_{Y_i}/\mathscr{I}_i^{d})\rightarrow
\text{H}^1(\PP^1,\mathscr{I}_i^d/\mathscr{I}_i^{d+1})$ be the
connection map and let
$\rho_{i,d}:\text{H}^0(\calO_{Y_{i}}/\mathscr{I}_{i}^{d+1})\rightarrow
\text{H}^0(\calO_{Y_{i}}/\mathscr{I}_{i}^{d})$. Consider the
following diagram
$$
\xymatrix@C=0.5cm{
  \text{H}^0(\calO_{Y_{i-1}}/\mathscr{I}_{i-1}^{d})
     &&\text{H}^0(\calO_{Y_{i}}/\mathscr{I}_{i}^{d})\ar[ll]_{\pi_{i,d-1}}
     &&\text{H}^0(\calO_{Y_{i+1}}/\mathscr{I}_{i+1}^{d})\ar[ll]_{\pi_{i+1,d-1}}\\
  \text{H}^0(\calO_{Y_{i-1}}/\mathscr{I}_i^{d+1})\ar[u]_{\rho_{i-1,d}}
     &&\text{H}^0(\calO_{Y_{i}}/\mathscr{I}_{i}^{d+1})\ar[ll]_{\pi_{i,d}}\ar[u]_{\rho_{i,d}}
     &&
}
$$
If both $\pi_{i,d-1}$ and $\pi_{i+1,d-1}$ are isomorphisms, then so
is $\pi_{i,d}$. To see this, note that surjectivity of
$\pi_{i+1,d-1}$ implies that $\delta_{i,d}=0$ by
Lemma\eqref{lifting}; and then we get that $\pi_{i,d}$ is
isomorphism by Lemma\eqref{iso}. Now, we can easily see that
$\pi_{i,0}:\text{H}^0(\PP^1,\calO_{Y_{i}}/\mathscr{I}_{i})\cong k
\rightarrow \text{H}^0(\PP^1,\calO_{Y_{i-1}}/\mathscr{I}_{i-1})\cong
k$ is isomorphism for all $i=1,2,\cdots$, hence all the
$\pi_{i,d}$'s are isomorphisms and all the $\delta_{i,d}$'s are 0.
This proves (a). To prove (b), we first note that
$\text{H}^0(\PP^1,\mathscr{I}/\mathscr{I}^2)$ is an $n-r$
dimensional $k$-vector space, say
$\text{Span}_k\{\alpha_1,\cdots,\alpha_{n-r}\}$, and it is a
subspace of $\text{H}^0(\PP^1,\calO_X/\mathscr{I}^2)$. Each
$\alpha_i$ lifts to a formal element $\hat{\alpha}_i\in
\text{H}^0(\mathfrak{X},\calO_\mathfrak{X})$ and this defines a ring
homomorphism
$$
\psi:k[[t_1,\cdots,t_{n-r}]]\rightarrow \text{H}^0(\mathfrak{X},\calO_\mathfrak{X}),
\quad t_i\mapsto \hat{\alpha}_i.
$$
To prove that $\psi$ is an isomorphism, we consider the following diagram
$$
\xymatrix{
  0\ar[r]
    &\text{H}^0(\mathscr{I}^d/\mathscr{I}^{d+1})\ar[r]
    &\text{H}^0(\calO_X/\mathscr{I}^{d+1})\ar[r]
    &\text{H}^0(\calO_X/\mathscr{I}^d)\ar[r]
    &0 \\
  0\ar[r]
    &\dfrac{(t_1,\cdots,t_{n-r})^d}{(t_1,\cdots,t_{n-r})^{d+1}}\ar[u]_{\lambda_d}\ar[r]
    &\dfrac{k[[t_1,\cdots,t_{n-r}]]}{(t_1,\cdots,t_{n-r})^{d+1}}\ar[u]_{\psi_{d+1}}\ar[r]
    &\dfrac{k[[t_1,\cdots,t_{n-r}]]}{(t_1,\cdots,t_{n-r})^d}\ar[u]_{\psi_d}\ar[r]
    &0
}
$$
and note that $\lambda_d$ factors through
$$
\dfrac{(t_1,\cdots,t_{n-r})^d}{(t_1,\cdots,t_{n-r})^{d+1}}\rightarrow
Sym^d(\text{H}^0(\mathscr{I}/\mathscr{I}^2))\rightarrow
\text{H}^0(Sym^d(\mathscr{I}/\mathscr{I}^2))\rightarrow
\text{H}^0(\mathscr{I}^d/\mathscr{I}^{d+1})
$$
and all these are isomorphisms and hence so is $\lambda_d$. Then we
prove that $\psi_d$ is isomorphism for all $d$ by induction since
$\psi_1$ is isomorphism. This shows
that $\psi$ is an isomorphism. $\square$\\

\textbf{Remark:} This theorem shows that there is a morphism of
formal schemes $\psi^a:\mathfrak{X}\rightarrow
\text{Spf}(k[[t_1,\ldots,t_{n-r}]])$. Then a natural question is:
Can we make this $\psi^a$ algebraic? The author expects that there
is an \'etale neighborhood $U$ of $\PP^1$ in $X$, i.e.
$$
\xymatrix{
 \PP^1\ar[r] &U\ar[r] &X
}
$$
such that, there is a morphism $U\rightarrow \mathbb{A}^{n-r}$ that
induces $\psi^a$. This is in the spirit of Artin's approximation
theorem, c.f.\cite{art} \cite{cdj}.

\section{Proof of Main Theorem}
In this section we prove our main theorem and give an application to
the Graber-Harris-Starr type theorem.
\begin{thm}(\textbf{Main Theorem})
Let $X/k$ be a quasi-projective variety over an algebraically closed
field $k$ of characteristic $p>0$. Assume that $X$ is freely
rationally connected. Then there is a separably rationally connected
variety $Y$ and a finite purely inseparable morphism $f:X\rightarrow
Y$. If $X$ is regular in codimension 1 (or normal), then so is $Y$.
\end{thm}

\textbf{Proof.} We want to produce the variety $Y$ by repeating
taking quotients by the foliation; By Lemma\eqref{stops}, we will
get a separably rationally connected variety after finitely many
steps if the ``splitting case" does not appear infinitely many
times. Then the regularity in codimension 1 (or normality) of $Y$
follows from Corollary \eqref{coreke}. So the proof of the main
theorem reduces to the following

\textbf{Claim: }\textit{Let $X$ be a separably uniruled variety over
$k$. If the quotient sequence \eqref{quot} has ``splitting case" in
each step,  then $X$ is not freely rationally connected.}

By contradiction, we assume that $X$ is FRC, i.e., there exists a
family of maximally free rational curves $\varphi:\PP^1\times
W\rightarrow X$ such that $\varphi^{(2)}:\PP^1\times\PP^1\times
W\rightarrow X^{(2)}=X\times X$ is dominant. First, we put the extra
assumption that each rational curve in this family is a closed
immersion. The fact that $\varphi^{(2)}$ being dominant implies that
there is a family of maximally free rational curves
$\varphi_0:\PP^1\times W_0\rightarrow X$ such that $\varphi_0$ is
dominant and $\varphi_0(0\times W_0)=x_0\in X(k)$. Actually, we can
just consider the following morphism $h$
$$
\xymatrix{
  h: \PP^1\times\text{PGL}_2\times W\ar[r] &\PP^1\times\PP^1\times W\ar[rr]^{\varphi^{(2)}} &&X^{(2)}
}
$$
with $h:(t,g,[\phi])\mapsto (\phi(g(0)),\phi(t))$. Take $\tilde{W}_0$ to be
$h^{-1}(x_0\times X)$.
$$
\xymatrix{
  \tilde{W}_0 \ar[r]\ar[d]^{h_0} &\PP^1\times\text{PGL}_2\times W\ar[d]^{h}\\
  x_0\times X \ar[r] &X\times X
}
$$
By definition, $h$ is dominant and hence $h_0$ is also dominant if we choose
$x_0\in X$ general enough. Then the following morphism is also dominant
$$
\tilde{h}_0:\PP^1\times \tilde{W}_0\to X,\quad (s,(t,g,[\phi]))\mapsto \phi(g(s)).
$$
By construction, we have $\tilde{h}_0(0\times\tilde{W}_0)=x_0$. Then we choose $W_0$
to be some component of $\tilde{W}_0$ and get a dominant
morphism $\varphi_0:\PP^1\times W_0\rightarrow X$ such that
$\varphi_0(0\times W_0)=x_0\in X(k)$. We denote $\PP^1\times W_0$ by
$Z$. Since $\varphi_0$ is dominant, there is a closed point
$z=(u,[\phi])\in Z$ such that $\calO_{X,x}\hookrightarrow
\calO_{Z,z}$ and hence $\hat{\calO}_{X,x}\hookrightarrow
\hat{\calO}_{Z,z}$, where $x=\varphi_0(z)$. Let
$\mathfrak{X}=X_{/\phi(\PP^1)}$ and
$\mathfrak{Z}=Z_{/\PP^1\times[\phi]}$. Then we have an induced
morphism of formal schemes $\psi:\mathfrak{Z}\to\mathfrak{X}$. By
Theorem\eqref{nbhd}, we have $\Gamma(\mathfrak{X},
\calO_\mathfrak{X})=k[[t_1,\cdots,t_{n-r}]]$ and it is easy to see
that $\Gamma(\mathfrak{Z},
\calO_\mathfrak{Z})=\hat{\calO}_{W_0,[\phi]}$. Now consider the
following diagram
\begin{equation}
\xymatrix{
  \Gamma(\mathfrak{X},\calO_\mathfrak{X}))\ar[rr]^{\psi^*}\ar[d]^{r_1} &&\Gamma(\mathfrak{Z},
     \calO_\mathfrak{Z})\ar[d]^{r_2} &&\\
  \hat{\calO}_{X,x} \ar[rr]^{\psi_z^*} &&\hat{\calO}_{Z,z}\ar[rr]^\cong
     &&\hat{\calO}_{\PP^1,u}\hat{\otimes}\hat{\calO}_{W_0,[\phi]}
}\label{injective}
\end{equation}
where the vertical maps are injective and $\psi_z^*$ is also
injective. So $\psi^*$ must be injective. But we can show that
$\psi^*(t_i)=0$. Indeed, the following factorization
$$
\xymatrix{
   &\Spec(k)\ar[dr]^{x_0} &\\
0\times W_0\ar[r]\ar[ur] &\PP^1\times W_0\ar[r] &X
}
$$
gives a factorization of $r_2\circ\psi^*$ as follows
$$
\xymatrix{
   &k\ar[dr] & & &\\
  \Gamma(\mathfrak{X},\calO_\mathfrak{X})\ar[ur]^{x_0^*}\ar[r] &\Gamma(\mathfrak{Z},\calO_\mathfrak{Z})\ar[r]
    &\hat{\calO}_{0\times W_0,(0,[\phi])}\ar[r]^\cong &\hat{\calO}_{W_0,[\phi]}\ar[r]
    &\hat{\calO}_{Z,z}
}
$$
In the above diagram $x_0^*(t_i)=0$ since $x_0\in \phi(\PP^1)$;
this gives a contradiction.

In the case where there is no maximally free rational curve on $X$
that is a closed immersion, we can find some integer $m>1$, such
that there is a maximally free rational curve $\phi=\phi_1\times
\phi_2 \times \ldots\times \phi_m$ on $X^{(m)}=X\times\cdots\times
X$ which is a closed immersion, where each $\phi_i$ is a maximally
free rational curve on $X$ for all $1\leq i\leq m$. It is also easy
to see that all maximally free rational curves on $X^{(m)}$ are of
the above form. Let
$\mathscr{D}^{(m)}:=p_1^*\mathscr{D}\oplus\cdots\oplus
p_m^*\mathscr{D}$. Since
$\text{Pos}(\phi^*(T_{X^{(m)}}))=\phi^*(\mathscr{D}^{(m)})$ for all
maximally free rational curves $\phi$ on $X^{(m)}$ and
$\mathscr{D}^{(m)}\subset T_{X^{(m)}}$ is saturated, we see that the
canonical foliation on $X^{(m)}$ is exactly $\mathscr{D}^{(m)}$. If
$Y$ is the quotient of $X$ by $\mathscr{D}$, then $Y^{(m)}$ is the
quotient of $X^{(m)}$. The quotient sequence of $X$ terminates with
an SRC variety if and only if the quotient sequence of $X^{(m)}$
does. Under the assumption that the quotient sequence of $X$ has
``splitting case" in each step, we get that the quotient sequence of
$X^{(m)}$ has ``splitting case" in each step. Thus we know that
$X^{(m)}$ is not FRC, hence $X$ is not FRC. This proves the the
claim and hence the main theorem.
$\square$\\

Now we are ready to prove the following application of the main
theorem.

\begin{thm}\label{fiberation}
Let $\pi: \mathscr{X}\to B$ be a proper flat family over a smooth
curve $B$, here everything is over an algebraically closed field $k$
of characteristic $p$. Assume that the geometric generic fiber of
$\mathscr{X}\to B$ is normal and freely rationally connected. Then
there is a morphism $s:B\to \mathscr{X}$ such that $\pi\circ
s=F_{\text{abs},B}^d$ for some $d\geq 0$, where
$F_{\text{abs},B}:B\to B$ is the absolute Frobenius morphism.
\end{thm}
\textbf{Proof.} We may shrink $B$ and assume that $B$ is affine. Let $\eta$ be the
generic point of $B$.

\textbf{Claim:} \textit{As a variety over $k$, $\mathscr{X}$ is separably uniruled}.\\
Indeed, let $\phi_{\bar{\eta}}:\PP^1_{\bar{\eta}}\to
\mathscr{X}_{\bar{\eta}}$ be a maximally free rational curve on the
geometric generic fiber. Assume that
$\phi_{\bar{\eta}}^*(T_{\mathscr{X}_{\bar{\eta}}})=\calO^{n-r}\bigoplus(\oplus
\calO(a_i))$ with $a_i\geq 1$ for $i=1,2,\ldots,r$. Then there is a
smooth curve $C/B$ such that $\phi_{\bar{\eta}}$ is actually defined
over $C$ and $\phi_C^*(T_{\mathscr{X}_C /C})$ splits uniformly as
$\calO^{n-r}\bigoplus(\oplus \calO(a_i))$. Namely, we then have the
following diagram
$$
\xymatrix{
  \PP^1_{\bar{\eta}}\ar[d]^{\phi_{\bar{\eta}}}\ar[r]
     &\PP^1_{C}\ar[d]^{\phi_{C}}
     &\\
  \mathscr{X}_{\bar{\eta}}\ar[r]\ar[d]
     &\mathscr{X}_C\ar[r]\ar[d]^{\pi_C}
     &\mathscr{X}\ar[d]^{\pi}\\
  \bar{\eta}\ar[r]
     &C\ar[r]
     &B
}
$$
By shrinking $C$, we may also assume that the image of $\phi_C$ is in
$(\mathscr{X}/B)^{sm}\times_B C$. Base change to a closed point
$c\in C(k)$, we get a free rational curve $\phi:\PP^1\to\mathscr{X}_b$ where
$b\in B(k)$ is the image of $c$. View $\phi$ as a rational curve on $\mathscr{X}$
and by construction we have the image of $\phi$ in
the smooth locus of $\mathscr{X}\to B$. Then we apply the following
short exact sequence
\begin{equation*}
\xymatrix{
  0\ar[r] &\phi^*(T_{\mathscr{X}_b})\ar[rr] &&
  \phi^*(T_{\mathscr{X}/k})\ar[rr] &&\phi^*\pi^*(T_{B/k})=\calO\ar[r] &0
}
\end{equation*}
Since $\phi^*(T_{\mathscr{X}_b})$ is globally generated,
$\phi^*(T_{\mathscr{X}/k})$ must be globally generated.

Let $\mathscr{D}$ be the canonical foliation on $\mathscr{X}$. By
construction, $\mathscr{D}$ is a submodule of $T_{\mathscr{X}/B}$.
Then $\mathscr{D}$ induces a foliation
$\mathscr{D}_{\bar{\eta}}\subset T_{\mathscr{X}_{\bar{\eta}}}$.
Since a free rational curve on $\mathscr{X}$ is unobstructed and moves to nearby fibres,
we have $\text{Pos}(\phi_{\bar{\eta}}^*T_{\mathscr{X}_{\bar{\eta}}})=\phi_{\bar{\eta}}^*\mathscr{D}_{\bar{\eta}}$
for a maximally free rational curve $\phi_{\bar{\eta}}$ on $\mathscr{X}_{\bar{\eta}}$. Since
$\mathscr{D}_{\bar{\eta}}\subset T_{\mathscr{X}_{\bar{\eta}}}$ is
saturated, we know that $\mathscr{D}_{\bar{\eta}}$ is the canonical
foliation on $\mathscr{X}_{\bar{\eta}}$.

By the Claim, we can construct the
quotient sequence of $\mathscr{X}$
\begin{equation}
\xymatrix{
  \mathscr{X}=\mathscr{Y}_0\ar[rr]^{f=f_0} &&\mathscr{Y}_1\ar[r]^{f_1}
  &\mathscr{Y}_2\ar[r]^{f_2} &\mathscr{Y}_3\ar[r]^{f_3}&\cdots
}\label{quot_X}
\end{equation}
Let $\mathscr{D}_i\hookrightarrow T_{\mathscr{Y}_i}$ be the
corresponding foliations. Note that all the $\mathscr{D}_i$ are in
vertical direction and hence $f_i$ are all defined over $B$. If we
base change the quotient sequence \eqref{quot_X} of $\mathscr{X}$ to
the geometric generic point $\bar{\eta}$ of $B$, we get exactly the
quotient sequence of $\mathscr{X}_{\bar{\eta}}$. Since
$\mathscr{X}_{\bar{\eta}}$ is FRC, our main theorem says that
$\mathscr{Y}_{i,\bar{\eta}}$ will eventually be an SRC normal
variety. Take $\mathscr{Y}=\mathscr{Y}_i$, where $i$ is large
enough. Then we have the following diagram
$$
\xymatrix{
  \mathscr{X}\ar[d]^{\pi}\ar[r]^f &\mathscr{Y}\ar[d]_{\tilde{\pi}}\ar[r]^g
     & \mathscr{X}^{(d)}\ar[r]^{\sigma^{(d)}}\ar[dl]^{\pi^{(d)}} &\mathscr{X}\ar[d]^{\pi}\\
  B\ar[r]^{=} &B\ar[rr]_{F_{\text{abs},B}^d} &&B
}
$$
where $\tilde{\pi}$ has normal SRC geometric generic fiber. By
\cite{js}, we can find a section of $\tilde{\pi}$, say $\tilde{s}$,
then we can just take $s=\sigma^{(d)}\circ g\circ\tilde{s}$.
$\square$

\begin{cor}
Let $X/k$ be a proper normal FRC variety over an algebraically
closed field $k$ of characteristic $p$. Then $X$ is simply
connected. Namely, the algebraic fundamental group $\pi_1(X)$ is
trivial.
\end{cor}
\textbf{Proof.} We have to show that every connected finite Galois
cover $\pi:Y\rightarrow X$ is trivial. Suppose that $\pi$ is
nontrivial. Since free rational curves on $X$ always lift to free
rational curves on $Y$, we know that $Y$ is also FRC. We can always
factor $\pi:Y\to X$ through $\pi':Y\to X'$ with $\pi'$ being cyclic
Galois cover with Galois group $G\cong \Z/\ell\,\Z$, where $\ell\in
\Z$ is a prime. To get a
contradiction, we only need to prove the following\\
\textit{Claim:} The action of $G$ on $Y$ has a fixed point.\\
To prove the claim, we fix an action of $G$ on $B'=\PP^1$ in the
following way:
\begin{itemize}
\item If $\ell\neq p$, then the generator $1\in G$ acts as $t\mapsto \zeta
t$ where $\zeta$ is a primitive $\ell^{\text{th}}$ root of unit.
\item If $\ell=p$, then the generator $1\in G$ acts as $t\mapsto
t+1$.
\end{itemize}
Then we have $B=B'/G\cong\PP^1$. Since $G$ acts on both $Y$ and
$B'$, we get a natural action of $G$ on $Y\times B'$ and let
$Z=Y\times B'/G$ be the quotient. Let $f:B'\rightarrow B$ and
$g:Y\times B'\rightarrow Z$ be the corresponding quotient morphisms.
Then we have the following commutative diagram
$$
\xymatrix{
  Y\times B'\ar[r]^{\quad g}\ar[d]_{p_2} &Z\ar[d]^{h}\\
  B'\ar[r]^{f} &B
}
$$
On the open part of $U\subset B$ where $G$ acts freely on
$f^{-1}(U)$, the above diagram is a fiber product square. Hence, for
a general $b\in B$, we have $Z_{b}=h^{-1}(b)\cong Y$ which is normal
and FRC. By Theorem \ref{fiberation}, we get $s:B\rightarrow Z$ such
that $h\circ s=F_{abs,B}^d$ is some power of the absolute Frobenius
morphism. Let $C=B'\times_{f,B,h\circ s}B$, and let
$\tilde{C}=\PP^1$ be the normalization of $C$. The the action of $G$
on $B$ induces an action on $C$ and hence also an action on
$\tilde{C}$. The morphism $s$ induces a $G$-equivariant morphism
$s':V\rightarrow Y\times B'$, where $V\subset C$ is the inverse
image of $U$. The morphism $s'$ induces a $G$-equivariant morphism
$\tilde{s}':\tilde{C}\rightarrow Y\times B'$. Then $\sigma=p_1\circ
\tilde{s}':\tilde{C}\rightarrow Y$ is also $G$-equivariant. Since
$\tilde{C}$ has at least one fixed point $x$, its image $\sigma(x)$
is a fixed point of $Y$. $\square$

\noindent Mingmin Shen\\
Department of Mathematics, Columbia University\\
mshen@math.columbia.edu

\end{document}